\documentclass[11pt, a4paper,oneside,reqno]{amsart}
\usepackage{amsfonts}
\usepackage{amssymb}
\usepackage{mathrsfs}
\usepackage{amsthm}
\usepackage{tcolorbox}
\usepackage{amsmath}
\usepackage{enumitem}
\usepackage{pifont}
\usepackage{float}
\usepackage{pgfplots}
\usepackage{mathtools}
\usepackage{cancel}
\usepackage{stmaryrd}
\usepackage{relsize}
\usepackage{xfrac}
\usepackage{anyfontsize}
\usepackage{color}
\usepackage{url}
\usepackage{bm}
\usepackage{datetime}
\usepackage{faktor}
\usepackage{tikz}
\usepackage[abs]{overpic}
\usepackage{tikz-cd}
\usepackage{forest}
\usetikzlibrary{calc, shapes, decorations.pathmorphing}
\usepackage{color}
\usepackage{scalerel}
\usepackage{graphicx}
\usepackage{csquotes}
\usepackage[top=1in, bottom=1in]{geometry}
\graphicspath{ {images/} }

\usepackage[colorlinks=false]{hyperref}
\usepackage{fancyhdr}
\allowdisplaybreaks
\usepackage[OT2,OT1]{fontenc}
\newcommand\cyr
{
\renewcommand\rmdefault{wncyr}
\renewcommand\sfdefault{wncyss}
\renewcommand\encodingdefault{OT2}
\normalfont
\selectfont
}
\DeclareTextFontCommand{\textcyr}{\cyr}

\pgfdeclarelayer{nodelayer}
\pgfdeclarelayer{edgelayer}
\pgfsetlayers{nodelayer,edgelayer,main}

\DeclarePairedDelimiter\ceil{\lceil}{\rceil}

\def \proof {\noindent\textit{Proof.}\hspace{5mm}}

\def \bs {$\hfill \blacksquare$\\}
\def \bsn {$\hfill \blacksquare$}
\def \sq {$\hfill \square$\\}
\def \sqn {$\hfill \square$}

\DeclareMathOperator{\Aut}{Aut}

\DeclareMathOperator{\Norm}{Norm}

\DeclareMathOperator{\ord}{ord}

\DeclareMathOperator{\Den}{Den}

\DeclareMathOperator{\ints}{INT}

\DeclareMathOperator{\Z}{\mathbb{Z}}
\DeclareMathOperator{\N}{\mathbb{N}}

\DeclareMathOperator{\p}{\mathfrak{p}}
\DeclareMathOperator{\q}{\mathfrak{q}}

\definecolor{supergreen}{RGB}{0, 170, 0}
\definecolor{superred}{RGB}{170, 0, 0}
\definecolor{darkred}{RGB}{100, 0, 0}
\definecolor{customcolour}{RGB}{0, 128, 128}

\newcommand{\Th}{\mbox{\ttfamily Th}}
\newcommand{\TTh}{\mbox{\emph{\ttfamily Th}}}
\newcommand{\ThE}{\mbox{\ttfamily Th}_{\exists}}
\newcommand{\ThAE}{\mbox{\ttfamily Th}_{\forall^1\exists}}
\newcommand{\TThE}{\mbox{\emph{\ttfamily Th}}_{\exists}}
\newcommand{\ThEp}{\mbox{\ttfamily Th}_{\exists^+}}
\newcommand{\TThEp}{\mbox{\emph{\ttfamily Th}}_{\exists^+}}

\newcommand{\fs}{\mathbb{F}_p [t]}
\newcommand{\fr}{\mathbb{F}_p (t)}

\newcommand{\iso}{\cong}

\newcommand{\defeq}{\mathrel{\mathop:}=}

\newtheorem{thm}{Theorem}[section]

\newtheorem{cor}[thm]{Corollary}
\newtheorem{lem}[thm]{Lemma}

\newtheoremstyle{case}{}{}{}{}{}{:}{ }{}
\theoremstyle{case}

\theoremstyle{definition}

\newtheorem{definition}[thm]{Definition}

\newtheorem{remark}[thm]{Remark}

\newtheorem{innercustomthm}{Theorem}
\newenvironment{customthm}[1]
  {\renewcommand\theinnercustomthm{#1}\innercustomthm}
  {\endinnercustomthm}
  
\newtheorem{innercustomcor}{Corollary}
\newenvironment{customcor}[1]
  {\renewcommand\theinnercustomcor{#1}\innercustomcor}
  {\endinnercustomcor}

\newtheorem{innercustomass}{Assumption}
\newenvironment{customass}[1]
  {\renewcommand\theinnercustomass{#1}\innercustomass}
  {\endinnercustomass}

\pgfplotsset{compat=1.14}

\makeatletter
\DeclareRobustCommand\widecheck[1]{{\mathpalette\@widecheck{#1}}}
\def\@widecheck#1#2{%
    \setbox\z@\hbox{\m@th$#1#2$}%
    \setbox\tw@\hbox{\m@th$#1%
       \widehat{%
          \vrule\@width\z@\@height\ht\z@
          \vrule\@height\z@\@width\wd\z@}$}%
    \dp\tw@-\ht\z@
    \@tempdima\ht\z@ \advance\@tempdima2\ht\tw@ \divide\@tempdima\thr@@
    \setbox\tw@\hbox{%
       \raise\@tempdima\hbox{\scalebox{1}[-1]{\lower\@tempdima\box
\tw@}}}%
    {\ooalign{\box\tw@ \cr \box\z@}}}
\makeatother

\begin{document}
\title{A Note on Hilbert's ``Geometric'' Tenth Problem}
\author{Brian Tyrrell}
\address{Mathematical Institute, Woodstock Road, Oxford OX2 6GG.}
\email{brian.tyrrell@maths.ox.ac.uk}

\begin{abstract}
This paper explores undecidability in theories of positive characteristic function fields in the ``geometric" language of rings $\mathcal{L}_F = \{0, 1, +, \cdot, F\}$, with a unary predicate $F$ for nonconstant elements. In particular we are motivated by a question of Fehm on the decidability of $\ThE(\fr; \mathcal{L}_F)$; equivalently, that of $\ThE(\fr; \mathcal{L}_r)$ \textit{without} parameters. We indicate how to generalise existing machinery to prove the undecidability of $\ThAE(K; \mathcal{L}_F)$ without parameters, where $K$ is the function field of a curve over an algebraic extension of $\mathbb{F}_p$, not algebraically closed. We discuss the problem (and its geometric implications) further in this context too.
\end{abstract}

\maketitle

\section{Introduction}
It appears to be a well known fact (e.g.\ \cite{denefpaper, pheidaszahidi}) that the diophantine problem for $\fs$ with coefficients in $\mathbb{F}_p$ is solvable. This is to say that Hilbert's Tenth Problem over $\fs$ with coefficients in $\mathbb{F}_p$ is solvable, or equivalently that $\ThEp(\fs)$ in the language of rings $\mathcal{L}_r$ \textit{without} parameters is decidable. The natural question to subsequently ask is whether Hilbert's Tenth Problem over the fraction field $\fr$ with coefficients in $\mathbb{F}_p$ is similarly solvable. This is more difficult to answer -- perhaps because in this context, the language of rings has an implicit geometric flavour.  

Indeed, consider the unary predicate $F$ defining the nonconstant elements of $\fr$, i.e.\ $\fr \models F(x) \iff x \in \fr\setminus\mathbb{F}_p$, and note that $F(\fr)$ is $\exists^+$-definable in the language of rings without parameters. Thus, we might equally have motivation to frame the above questions in the ``geometric'' language $\{0, 1, +, \cdot, F\}$. From the perspective of algebraic geometry, the subtlety between these problems is highlighted:

\begin{thm}\label{1point1}
There \emph{exists} an algorithm which upon input an $\mathbb{F}_p$-variety $\mathscr{V}$ outputs ``YES'' if there exists an $\mathbb{F}_p$-morphism $\mathbb{P}^1 \rightarrow \mathscr{V}$, and ``NO'' otherwise.

There \emph{does not exist} an algorithm which upon input an $\mathbb{F}_p$-variety $\mathscr{V}$ outputs ``YES'' if there exists a \emph{nonconstant} $\mathbb{F}_p$-morphism $\mathbb{P}^1 \rightarrow \mathscr{V}$, and ``NO'' otherwise.
\end{thm}

\proof
The former follows from the aforementioned decidability of $\ThEp(\fs; \mathcal{L}_r)$ without parameters. The latter is \cite[Theorem 2.1]{pheidaszahidi2}. \bs

In this paper, we will seek to answer decidability questions in the geometric language. In $\S \ref{machinery}$ we generalise techniques of Eisentr\"{a}ger \&\ Schlapentokh and use modern results of Pasten (recalled in $\S \ref{recap}$) to conclude undecidability in the language of rings augmented by a predicate $B_{l}$ defining ``good behaviour'' relative to a fixed prime $l$ (elements that are \emph{$l$-behaved}; see \emph{Definition \ref{lfriend}}). We are unable to prove \textit{existential} undecidability results ($\S \ref{break}$ explains why) however we are able to get some distance by using a single initial universal quantifier:

\begin{customthm}{\ref{firstuni}}
\textit{Let $K$ be the function field of a curve, of positive characteristic and not containing the algebraic closure of a finite field. Then for some prime $l$, $\TTh_{\forall^1\exists^+}(K; \mathcal{L}_F \cup \{\neg B_{l}\})$ without parameters is undecidable.}
\end{customthm}

We then make some progress in $\S \ref{universal}$ eliminating this introduced behaviour-controlling predicate $B_l$, to obtain as a corollary:

\begin{customcor}{\ref{firstunicor}}
\textit{Let $K$ be the function field of a curve, of positive characteristic and constant subfield $C \subsetneq \widetilde{\mathbb{F}_p}$. Then $\TTh_{\forall^1\exists}(K; \mathcal{L}_F)$ without parameters is undecidable.}
\end{customcor}

\noindent As in \emph{Theorem \ref{1point1}}, this can be reformulated in geometric terms.

\begin{customcor}{\ref{firstunicor} (restated)}
\textit{Let $\mathcal{C}$ be a curve defined over a field $\mathbb{F}$, $\mathbb{F}_p \subseteq \mathbb{F} \subsetneq \widetilde{\mathbb{F}_p}$ for some prime $p$. Given an $\mathbb{F}$-morphism $\pi: \mathscr{V} \rightarrow \mathscr{W}$ of $\mathbb{F}$-varieties with $\mathscr{W} \subseteq \mathbb{P}^1$ infinite, it is undecidable when one can prove for all nonconstant $\mathbb{F}$-rational maps $r : \mathcal{C} \dashrightarrow \mathscr{W}$ there exists an $\mathbb{F}$-rational map $\mathcal{C} \dashrightarrow \mathscr{V}$ making the below diagram commute.}
\[
\begin{tikzcd}
& \mathcal{C} \arrow[ld, dashed, "\exists", swap] \arrow[rd, dashed, "r"] & \\
\mathscr{V} \arrow[rr, "\pi"] \arrow[rr, phantom, ""{name=D, above}]{}& & \mathscr{W} \arrow[from=1-2, to=D, "\circlearrowleft", phantom]
\end{tikzcd}
\]
\end{customcor}

Following a suggestion due to E.\ Hrushovski, in $\S \ref{concl}$ we demonstrate a scenario where a natural subset of $l$-behaved elements is existentially $\mathcal{L}_F$-definable without parameters (this will be when the function field has sufficiently high genus, indicating the use of the predicate $B_l$ is in a sense ``not too strong'' and only applicable to low genus function fields). Finally we demonstrate how to use the machinery established here to (re)prove\footnote{The author suspects \emph{Corollary \ref{three12}} (in model-theoretic terms: \emph{Theorem \ref{rephrasing}}) is known, but is not aware of a source for it.} the following geometric result:

\begin{customcor}{\ref{three12}}
\textit{Let $\mathcal{C}$ be a curve defined over a finite field $\mathbb{F}_p$ \emph{(}$p > 2$\emph{)} with a nonconstant $\mathbb{F}_p$-rational map $\tau : \mathcal{C} \dashrightarrow \mathbb{P}^1$. No algorithm exists which, upon input an $\mathbb{F}_p$-variety $\mathscr{A}$ and an $\mathbb{F}_p$-rational map $\rho : \mathscr{A} \dashrightarrow \mathbb{P}^1$, outputs ``YES'' if there exists an $\mathbb{F}_p$-rational map $\mathcal{C} \dashrightarrow \mathscr{A}$ making the below diagram commute, and ``NO'' otherwise. }
\[
\begin{tikzcd}
 \mathcal{C} \arrow[rr, dashed, "\exists"] \arrow[rd, dashed, "\tau", swap]& \arrow[d, phantom, "\circlearrowleft"] & \mathscr{A} \arrow[ld, dashed, "\rho"]\\
 & \mathbb{P}^1 &
\end{tikzcd}
\]
\end{customcor}

\subsection{Interlude: rewriting equations over finite extensions}
We will require the following standard result on diophantine sets in field extensions, found in \cite{eisenshlap}.

\begin{thm}\label{ugly}
Let $K$ be a field and let $f(T_1, \dots, T_n, X_1, \dots X_{n_2}, Y_1, \dots, Y_{n_3})$, \linebreak $g(X, T_1, \dots, T_n)$ be polynomials with coefficients in $K$. Assume the degree of $g$ in $X$ is positive and the same for all values of $T_1, \dots, T_n$ (i.e.\ the leading coefficient of $g$ as a polynomial in $X$ over the algebraic closure of $\mathbb{F}_p(T_1, \dots, T_n)$ is always nonzero for any choice of $T_1, \dots, T_n \in K$, and the degree of $g$ in $X$ is positive). Let $A \subseteq K^n$ be defined as follows:
\begin{align*}
(t_1, \dots t_n) \in A \iff \exists& x_1, \dots x_{n_2} \in K\mbox{, } x \in \widetilde{K}\mbox{, } y_1, \dots, y_{n_3} \in K(x) \mbox{ s.t. }\\
&g(x, t_1, \dots, t_n) = 0 \land f(t_1, \dots, t_n, x_1, \dots, x_{n_2}, y_1, \dots, y_{n_3}) = 0.
\end{align*}

Then $A$ has a diophantine definition over $K$. Moreover, there is a diophantine definition of $A$ with coefficients depending only on the coefficients and degrees of $g$ and $f$, and this definition can be constructed effectively from these coefficients.
\end{thm}

\proof
See \cite[Lemma B.7.5]{shlapbook} or \cite[Lemma 1.3]{shlapen}. \bsn

\bigskip
\section{Machinery}\label{machinery}
In modern publications (cf.\ \cite{eisenshlap, shlap2, shlapen, shlapbook}) there is a standard two step process to conclude undecidability of the function field $K$ of a curve $\mathcal{C}$ with field of constants $C$ in characteristic $p > 0$. Denote $C_0 \defeq \widetilde{\mathbb{F}_p} \cap C$.

\begin{enumerate}
    \item Define the ``Denef'' predicate $\Den_p(x, y) \iff \exists s \in \N \mbox{ } x = y^{p^s} \lor y = x^{p^s} \iff \exists s \in \Z \mbox{ } x = y^{p^s}$.\label{steps}
    \item Define for a nonconstant `special' element $u \in K$ and any prime $\p$ with $\ord_{\p} u = 1$, the set $\ints(K, \p, u)$, where if $x \in \ints(K, \p, u)$ then $\ord_{\p} x \geq 0$, and if $x \in C_0 (u)$ and $\ord_{\p} x \geq 0$ then $x \in \ints(K, \p, u)$.
\end{enumerate}

\noindent This Denef predicate is useful by the following theorem:

\begin{thm}\label{modpheidas}
Let $p$ be a prime number. Then $\TThEp(\N; 0, 1, +, |_p, \leq)$ is undecidable, where $|_p$ is the binary predicate defined by $a |_p b \iff \exists s \in \N, \mbox{ }b = p^s a \lor a = p^s b$.
\end{thm}

\proof
Pheidas in \cite{pheidas87} demonstrated $\ThEp(\N; 0,1, +, |^p)$ is undecidable, where now \linebreak $n |^p m \iff \exists s \in \N, \mbox{ } m = p^s n$. This relation is definable in $\ThEp(\N; 0,1, +, |_p, \leq)$, as 
$$x |^p y \iff x|_p y \land x \leq y.$$
We conclude $\ThEp(\N; 0, 1, +, |_p, \leq)$ is undecidable, as required. \bs

Therefore if one can interpret the theory $\ThEp(\N;0,1,+, |_p, \leq)$ in $\ThE(K)$, one can conclude undecidability.

\subsection{Definability of the Denef Predicate.}\label{recap}
Step ($1$) is concluded in our setting (in the language $\mathcal{L}_F$, no parameters from the underlying field permitted) using arithmetic properties detailed by Pasten \cite{pastennew} for all odd primes $p$ (in fact, Pasten's method produces an $\mathcal{L}_F$-formula for the Denef predicate \textit{uniform} in such $p$). We begin with the below theorem -- all results in this subsection are obtained from Pasten \cite{pastennew}, and elaborated here for the sake of exposition, that parameters from $K$ are not necessary with $\mathcal{L}_F$.

\begin{thm}\label{above}
Let $g \geq 0$, $d \geq 1$ be integers and let $p > 2$ be a prime. Let
$$M = M(g, d, p) = \left\lceil\tfrac{1}{d}\left(4g + 12 + 8 \sum_{i = 1}^{\ceil{(d-1)/2}} p^i \right)\right\rceil.$$
Then we have the following:

Let $k$ be a field of characteristic $p$ and let $K$ be a one variable function field of genus $g$ defined over $k$. Let $F_1, \dots, F_M \in \mathbb{F}_p[X]$ be pairwise coprime irreducible polynomials of degree $d$. Take $f, h \in K$ both nonconstant. There exists $s \in \N$ such that $f = h^{p^s}$ or $h = f^{p^s}$ if and only if $F_i(f) F_i(h)$ is a square in $K$ for each $i = 1, \dots, M$.
\end{thm}

\proof
The forward implication is \cite[Remark 1]{pastennew}, while the reverse implication is \cite[Theorem 1.6]{pastennew} exactly. \bs

\begin{thm}
Let $g \geq 0$ be an integer. There exists an existential $\mathcal{L}_F$-formula $\varphi_g(x, y)$ with the following property:

Given any prime $p > 2$, any field $k$ of positive characteristic $p$, and any function field $K/k$ of a curve of genus $\leq g$, for every pair of elements $f, h \in K$,
$$K \models \varphi_g(f, h) \iff \exists s \in \N \mbox{ s.t. } f = h^{p^s} \mbox{ or } h = f^{p^s}.$$
\end{thm}

\proof
This is a minor adaption of \cite[Theorem 1.5]{pastennew} to remove the dependency on the parameter. For $d \geq 1$, let $M_d$ be the number of monic irreducible polynomials of degree $d$ in $\mathbb{F}_p[X]$. One may first prove for sufficiently large\footnote{Pasten notes that $d \geq 2 \log (12 + \sqrt{8g + 168})$ suffices.} primes $d$, $M_d > M = M(g, d, p)$ from \emph{Theorem \ref{above}}. Therefore we may always chose distinct monic irreducible polynomials $F_1, \dots, F_M \in \mathbb{F}_p[X]$ as required for \emph{Theorem \ref{above}}.

Let $\phi_{g, p}(x, y)$ be the formula
$$\bigwedge_{i = 1}^M \exists z \mbox{ }(\widetilde{F}_i(x) \widetilde{F}_i(y) = z^2),$$

where $\widetilde{F}_i$ is a lift of $F_i$ from $\mathbb{F}_p[X]$ to $\Z[X]$. Let $k'$ be the constant subfield of $K$. By \emph{Theorem \ref{above}}, if $K \models \phi_{g,p}(f, h) \lor (\neg F(f) \land \neg F(h))$ then either $f, h \in k'$ or $\exists s \in \N$ s.t.\ $f = h^{p^s}$ or $h = f^{p^s}$. Conversely, if $\exists s \in \N$ s.t.\ $f = h^{p^s}$ or $h = f^{p^s}$, then $K \models \phi_{g,p}(f, h)$.

Note that if $p > 4g + 12$ we may take $d = 1$ (hence $M = 4g + 12$) and chose polynomials $F_i = X - i$ for $i = 1,\dots, M$. Hence for $p > 4g +12$, we may chose $\phi_{g,p}$ uniformly in $p$; denote this $\phi_g$. 

Let $\chi_g(x, y)$ be the formula

$$\left(\left[\phi_g(x, y) \land \bigwedge_{p \leq 4g + 12} p \neq 0\right] \lor \left[\bigvee_{p \leq 4g + 12} \{p = 0 \land \phi_{g, p}(x, y)\}\right]\right),$$

\noindent and $\varphi_g(x, y)$ be the formula
\begin{align*}
\big( F(x) \land F(y) \land \chi_g(x, y) \big) \lor \big( \neg F(x) \land \neg F(y) \land \exists u, v \big[ F(u) \land F(v) \land \chi_g(u, v) \land \chi_g(ux, vy)\big]\big)
\end{align*}

We claim this is the required formula. Indeed, if $f = h^{p^s}$ or $h = f^{p^s} \in K$ for $s \in \N$, then either the first or second disjunct of $\varphi_g$ is satisfied, depending on whether $f, h$ are both constant or not. On the other hand, if $K \models \varphi_g(f, h)$, then either $f, h$ are nonconstant and $f = h^{p^s}$ or $h = f^{p^s}$ for $s \in \N$, or $f, h$ are constant and there exists $s_1, s_2 \in \N$, $u, v \in K \setminus k'$ such that $u = v^{p^{s_1}}$ or $v = u^{p^{s_1}}$, and $uf = (vh)^{p^{s_2}}$ or $vh = (uf)^{p^{s_2}}$.

Hence, either $f = v^{p^{s_2} - p^{s_1}} h^{p^{s_2}}$, or $h = u^{p^{s_2} - p^{s_1}} f^{p^{s_2}}$, or $h = v^{p^{s_1+s_2}-1} f^{p^{s_2}}$, or $f = u^{p^{s_1+s_2}-1} h^{p^{s_2}}$. In the former two cases, if $s_1 \neq s_2$, this forces $v$ or $u$ to be constant; a contradiction. Therefore either $f = h^{p^{s_2}}$ or $h = f^{p^{s_2}}$ in this case. 

In the latter two cases, if either $s_1$ or $s_2$ is not $0$, this again forces $v$ or $u$ to be constant; a contradiction. Therefore in this case $f = h = h^{p^{0}}$, as desired. \bs

\begin{cor}\label{pastenbuchi}
Let $K$ be the function field of a curve of genus $g$ over a field of characteristic $p >2$. The Denef predicate $\Den_p$ is existentially $\mathcal{L}_F$-definable in $K$, without parameters. In particular, Step \emph{(1)} of p.\! \pageref{steps} is immediately concluded. \bsn
\end{cor}

\bigskip
\subsection{Integrality.}\label{integrality}

The ideas behind the machinery presented here originate with Er{\v s}ov \cite{ersov} and Penzin \cite{pen}, J.\ Robinson \cite{robinson, rob59} and Rumely \cite{rumely}. They have been used extensively in this context by Shlapentokh (e.g.\ \cite{shlapen1, shlap2, shlap103, shlap101, shlap4, shlap3, shlap}, and explained in detail in \cite[Chapter 4]{shlapbook}), and Eisentr\"{a}ger-Shlapentokh \cite{eisenshlap}. Indeed, the process below is identical to \cite[\S 6]{eisenshlap} with the exception that here we draw attention to any parameters used. The goal is to produce an existentially $\mathcal{L}_r$-defined subset of any function field, which forces its members to have ``small'' poles. We will use a parameter in defining this set, but keep careful control on what exact properties this parameter satisfies.

Let $p$ be any prime, and $K$ be a (one variable) function field of characteristic $p > 0$ with any subfield of constants $C$. Let $l$ be a prime, not necessarily distinct to $p$. 

\begin{definition}\label{lfriend}
An element $u \in K$ is \textit{$l$-behaved} if there exists a prime $\p$ of $K$ such that $v_{\p}(u) > 0$, $v_{\p}(u) \not\equiv 0 \mod l$, and $[Kv_{\p} : C] \not\equiv 0 \mod l$. Define:
$$\mathfrak{z}_b(u) \defeq \prod \{ \p^{v_{\p}(u)} \mbox{ : }  v_{\p}(u) > 0\mbox{, } v_{\p}(u) \not\equiv 0\!\! \mod l\mbox{, }[Kv_{\p} : C] \not\equiv 0\!\! \mod l\},$$
the \textit{$l$-behaved factor} of the zero divisor $\mathfrak{z}(u)$ of $u$.
\end{definition}

\noindent To begin, we assume the following:

\begin{customass}{({\Large$\star$})}\label{assump}
$C_0$, the algebraic closure of $\mathbb{F}_p$ in $C$, contains an $l$-th primitive root of unity for some prime $l \neq p$, and for some $a \in C_0$, $K$ will not contain any root of $T^l - a$.
\end{customass}

E.g.\ if $K$ is global, $C_0$ will be a finite field $\mathbb{F}_{p^n}$; ({\Large$\star$}) amounts to assuming $l | (p^n - 1)$. Let $u \in K$ be $l$-behaved. In addition, fix the following notation:
\begin{itemize}
\item Let $K_0$ denote the algebraic closure of $C_0(u)$ in $K$.
    \item Let $\delta \in \widetilde{K}$ be a root of $T^l -(u+1)$.
    \item For $w \in K$, let $h_w = \tfrac{w^l}{u} + \tfrac{1}{u^l} = \tfrac{w^l u^{l-1} + 1}{u^l}$.
    \item Let $\beta_w \in \widetilde{K}$ be a root of $T^l - (\tfrac{1}{h_w} +1)$.
    \item Let $\alpha \in \widetilde{K}$ be a root of $T^l - a$.
\end{itemize}

\noindent Consider the following extensions:

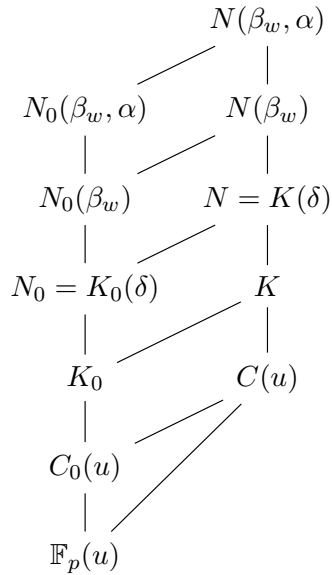
\begin{center}
\begin{tikzcd}[column sep=tiny, row sep=scriptsize]
& N(\beta_w, \alpha) \arrow[d, dash]\\
N_0(\beta_w, \alpha) \arrow[d, dash] \arrow[ru, dash] & N(\beta_w) \arrow[d, dash]\\
N_0(\beta_w) \arrow[d, dash] \arrow[ru, dash] & N = K(\delta) \arrow[d, dash]\\
N_0 = K_0(\delta) \arrow[d, dash] \arrow[ru, dash]& K \arrow[d, dash]\\
K_0 \arrow[ru, dash] \arrow[d, dash] &C(u) \arrow[dl, dash]\\
C_0(u)  \arrow[d, dash] & \\
\mathbb{F}_p(u) \arrow[ruu, dash]&\\
\end{tikzcd}    
\end{center}

We have the following series of lemmas from \cite{eisenshlap}. From this point onward, we will interchange the notations ``$v_{\p}$'' and ``$\ord_{\p}$'' for a prime $\p$.

\begin{lem}\label{wellillbedamned}
Let $G$ be a field of positive characteristic $p$ possessing a primitive $l$-th root of unity $\xi_l$. Let $\alpha \in \widetilde{G}$, let $\alpha_j = \xi_l^j\alpha$, $j = 0, \dots , l - 1$, and let
$$P(a_0, \dots, a_{l-1}) = \prod_{j=0}^{l-1} (a_0 + a_1 \alpha_j + \dots + a_{l-1} \alpha^{l-1}_j).$$

In this case, if $[G(\alpha) : G] = l$, then $P(a_0, \dots , a_{l-1}) = \Norm_{G(\alpha)/G} (a_0 + a_1 \alpha + \dots + a_{l-1} \alpha^{l-1})$. Also, if $\alpha \in G$, then for any $y \in G$ the equation $P(X_0, \dots , X_{l-1}) = y$ has solutions $x_0, \dots, x_{l-1} \in G$. \bsn
\end{lem}

\begin{lem}\label{mention}
Let $G/H$ be a Galois extension of algebraic function fields of degree $k$. Let $\p$ be a prime of $H$ with only one, unramified, factor in $G$. Let $x \in H$ be such that $\ord_{\p}x \not\equiv 0 \mod k$. Then $x$ is not a norm of an element of $G$. \bsn
\end{lem}

\begin{lem}\label{locallem}
Let $G/H$ be an unramified extension of local fields of degree $k$. Let $\mathfrak{m}$ be the prime of $H$. Let $x \in H$ be such that $\ord_{\mathfrak{m}} x \equiv 0 \mod k$. Then $x$ is a norm of some element of $G$.\bsn
\end{lem}

\begin{lem}\label{mainlem}
Let $L$ be a function field of characteristic $p$ possessing an $l$-th primitive root of unity. Let $z \in L$ and let $\gamma$ be a root of the equation $X^l - z = 0$. If for some $L$-prime $\mathfrak{a}$, $\ord_{\mathfrak{a}}(z) < 0$ and $\ord_{\mathfrak{a}}(z) \not\equiv 0$ modulo $l$ then $\mathfrak{a}$ is completely ramified in $L(\gamma)/L$. Also, if $z$ is integral at $\mathfrak{a}$ and $z$ is equivalent to a nonzero $l$-th power modulo $\mathfrak{a}$, then $\mathfrak{a}$ will split completely in $L(\gamma)$. \bs
\end{lem}

Applied to the fields in question, they have the following corollary:

\begin{cor}\label{thecor}
The following statements are true about the extensions $N/K$, $N_0/K_0$:
\begin{enumerate}
    \item There is no constant field extension.
    \item The factors of $\mathfrak{z}_b(u)$ split completely, into factors of relative inertial degree 1.
    \item Any factor $\q$ of $\mathfrak{pl}(u)$ (the pole divisor of $u$) where $\ord_{\q} (u) \not\equiv 0 \mod l$ ramifies completely into factors of inertial degree 1.
\end{enumerate}
In particular, when $[N:K] > 1$ $($resp.\ $[N_0 : K_0] > 1)$, for any prime $\widehat{\q}$ of $N$ $($resp.\ $N_0)$ over $\mathfrak{pl}(u)$, $\ord_{\widehat{\q}}(u) \equiv 0 \mod l$. \bs
\end{cor}

\begin{lem}\label{thelem}
For $w \in K$, the following statements are true in $N(\beta_w)$ $($for $w \in K_0$, we obtain the analogous statements for $N_0(\beta_w)/N_0/K_0)$:
\begin{enumerate}
    \item If $\widehat{\p}$ is a prime of $N(\beta_w)$ and $\widehat{\p} | \mathfrak{z}_b(u)$ in $N(\beta_w)$ while $\ord_{\widehat{\p}}(w) < \tfrac{1-l}{l}\ord_{\widehat{\p}}(u)$, then \linebreak $f(\widehat{\p}/\overline{\p}) = f(\overline{\p}/\p) = 1$, where $\overline{\p} = \widehat{\p} \cap N$ and $\p = \overline{\p} \cap K$.
    \item If $\widehat{\p}$ is a prime of $N(\beta_w)$ and $\widehat{\p} | \mathfrak{z}_b(u)$ in $N(\beta_w)$ while $\ord_{\widehat{\p}}(w) < \tfrac{1-l}{l}\ord_{\widehat{\p}}(u)$, then \linebreak $\ord_{\widehat{\p}}(h_w) \not\equiv 0 \mod l$.
    \item If $\mathfrak{t}$ is a prime of $N(\beta_w)$ and $\mathfrak{t} \nmid \mathfrak{z}(u)$, then $\ord_{\mathfrak{t}}(h_w) \equiv 0 \mod l$. 
    \item If $\p$ is a prime of $K$ such that $\p | \mathfrak{z}_b(u)$ and $\ord_{\p}(w) \geq \tfrac{1-l}{l}\ord_{\p}(u)$, then $\ord_{\p}(h_w) \equiv 0 \mod l$.
\end{enumerate}
\end{lem}

\proof 
First, a calculation. Suppose $\p$ is a prime of $K$ such that $\p | \mathfrak{z}_b(u)$. Then:
\begin{align*}
\ord_{\p}(h_w) &= \ord_{\p}(w^l u^{l-1} + 1) - l\ord_{\p}(u)\\
&= \min\{l\ord_{\p}(w) + (l-1)\ord_{\p}(u), 0\} - l\ord_{\p}(u),
\end{align*}

as $\ord_{\p}(u) \not\equiv 0 \mod l$. If $\ord_{\p}(w) < \tfrac{1-l}{l}\ord_{\p}(u)$, then $\ord_{\p}(w^lu^{l-1} + 1) = l\ord_{\p}(w) + (l-1)\ord_{\p}(u)$, and $\ord_{\p}(h_w) \not\equiv 0 \mod l$.

Now, (1) \&\ (2). Let $\p | \mathfrak{z}_b(u)$ in $K$ and note by \emph{Corollary \ref{thecor}}, $\p$ splits completely in $N$ into factors of inertial degree 1. For any primes $\widehat{\p} | \overline{\p} | \p$ in $N(\beta_w)$ (resp.\ $N$), we have $\ord_{\widehat{\p}}(h_w) = \ord_{\overline{\p}}(h_w) = \ord_{\p}(h_w)$ and $f(\widehat{\p}/\overline{\p}) = f(\overline{\p}/\p) = 1$, as by \emph{Lemma \ref{mainlem}}, $\p$ and $\overline{\p}$ split completely in their extensions. We reach the desired conclusion for (2) from our initial calculation.

Part (3) is \cite[Lemma 6.9(3)]{eisenshlap} exactly. Part (4) is a straightforward calculation based on the first paragraph. \bs

We now replicate the sufficient and necessary conditions of \cite{eisenshlap} for $w \in K$ to have ``small'' poles at factors $\p$ of $\mathfrak{z}_b(u)$, in terms of a norm equation. 

\begin{lem}\label{thing1}
If $\ord_{\p}(w) < \tfrac{1-l}{l}\ord_{\p}(u)$ at any factor $\p | \mathfrak{z}_b(u)$ in $K$, then there is no $x \in N(\beta_w, \alpha)$ such that $\Norm_{N(\beta_w, \alpha)/N(\beta_w)}(x) = h_w$.
\end{lem}

\proof
Let $\widehat{\p}$ be a factor of $\mathfrak{z}_b(u)$ in $N(\beta_w)$ over $\p$; then $\ord_{\widehat{\p}}(w) = \ord_{\p}(w) < \tfrac{1-l}{l}\ord_{\p}(u)$ by the argument of \emph{Lemma \ref{thelem} (1)}. By the same reasoning as \emph{Corollary \ref{thecor} (1)}, there is no constant field extension in $N(\beta_w)/N$, and since $f(\widehat{\p}/\p) = 1$, the equation 
\begin{equation}\label{theeq}
T^l - a = 0    
\end{equation}

has no root in the residue field of $\widehat{\p}$ in $N(\beta_w)$ if and only if (\ref{theeq}) has no root in the residue field of $\p$ in $K$. Indeed, this is the case, as by design (\ref{theeq}) has no root in $K$ (hence $C$), and  $[Kv_{\p} : C] \not\equiv 0 \mod l$. Therefore $\widehat{\p}$ cannot split in the extension $N(\beta_w, \alpha)/N(\beta_w)$, as the extension is of prime degree and the residue field of $\widehat{\p}$ must extend. If $h_w$ is to be a norm in this extension, then $\ord_{\widehat{\p}}(h_w) \equiv 0 \mod l$ by \emph{Lemma \ref{mention};} however by \emph{Lemma \ref{thelem} (2)}, we see that $\ord_{\mathfrak{p}}(h_w) \not\equiv 0 \mod l$, and by the argument of \emph{Lemma \ref{thelem}}, $\ord_{\widehat{\mathfrak{p}}}(h_w) = \ord_{\mathfrak{p}}(h_w)$. This is a contradiction, as desired. \bs

\begin{lem}\label{thing2}
For $w \in C_0(u)$, if $\ord_{\p}(w) \geq \tfrac{1-l}{l}\ord_{\p}(u)$ for all factors $\p | \mathfrak{z}_b(u)$ in $K$, then there exists $x \in N(\beta_w, \alpha)$ such that $\Norm_{N(\beta_w, \alpha)/N(\beta_w)}(x) = h_w$.
\end{lem}

\proof
First, note that it is sufficient to prove there exists $x \in N_0(\beta_w, \alpha)$ such that $\Norm_{N_0(\beta_w, \alpha)/N_0(\beta_w)}(x) = h_w$. Indeed, $x \in N_0(\beta_w, \alpha)$ has the same $N_0(\beta_w)$-coordinates with respect to the power basis of $\alpha$ as in $N(\beta_w)$, hence $x$ has the same conjugates over $N(\beta_w)$ and $N_0(\beta_w)$, hence $x$ has the same norm.

Next we claim the divisor of $h_w$ is an $l$-th power of another divisor in $N_0(\beta_w)$. We may assume $N_0(\beta_w)/N_0$ is a proper extension; otherwise $\tfrac{h_w+1}{h_w}$ is an $l$-th power in $N_0(\beta_w) = N_0$ and we immediately conclude the claim. Whenever there is a prime $\q$ such that $\ord_{\q}(u) \leq 0$, the claim follows from a simple calculation based on the construction of $h_w$ and using \emph{Corollary \ref{thecor}}. Note the assumption that $w \in C_0(u)$ forces $\ord_{\p}(w) = e(\p/u)\ord_u(w) = \ord_{\p}(u)\ord_u(w)$. We are also assuming for $\p | \mathfrak{z}_b(u)$ we have $\ord_{\p}(w) \geq \tfrac{1-l}{l}\ord_{\p}(u)$. 

By these points above, it cannot be the case that for some prime $\q$, $\ord_{\q}(u) > 0$ and $\ord_{\q}(w) < 0$: otherwise $\ord_{\p}(w) \leq -\ord_{\p}(u)$ for all $\p | \mathfrak{z}_b(u)$, hence
$$-(l-1)\ord_{\p}(u) < l \ord_{\p}(w) \leq -l \ord_{\p}(u),$$
a contradiction.

If $\ord_{\q}(u) > 0$, $\ord_{\q}(w) = 0$, then by calculation $\ord_{\q}(h_w) \equiv 0 \mod l$.

Finally, if for some prime $\q$ it is the case that $\ord_{\q}(u) > 0$, $\ord_{\q}(w) > 0$, then $\ord_{\q}(w) = \ord_{\q}(u) \ord_u(w) > 0$, hence $\ord_u(w) > 0$, thus $\ord_{\p}(w) \geq \ord_{\p}(u)$ for all $\p | \mathfrak{z}_b(u)$, and therefore $\ord_{\p}(h_w) \equiv 0 \mod l$ as required.

Observe there is a \emph{finite} extension $\widehat{N}_0/\mathbb{F}_p(u)$ such that $w \in \widehat{N}_0$, the divisor of $h_w$ is an $l$-th power of another divisor of $\widehat{N}_0$, and $\alpha$ is of degree $l$ over $\widehat{N}_0$. By a similar argument as in the beginning of the lemma, it is sufficient to solve $\Norm_{\widehat{N}_0(\beta_w, \alpha)/\widehat{N}_0(\beta_w)}(x) = h_w$. By \emph{Lemma \ref{mainlem}}, the extension $\widehat{N}_0(\beta_w, \alpha)/\widehat{N}_0(\beta_w)$ is unramified, hence by Weil \cite[Corollary, p.\ 226]{weil} \textit{locally} every unit is a norm.  Therefore, by \emph{Lemma \ref{locallem}} and the Strong Hasse Principle\footnote{For this we require our fields to be global.} \cite[Theorem 32.9]{reiner}, \emph{globally} $h_w$ is a norm, as required. \bs

\noindent We finish this argument with the following theorem and corollary:

\begin{thm}\label{thethm}
Let $\alpha_j = \xi_l^j\alpha$ for $j = 0, \dots, l-1$, and let
$$P(a_0, \dots, a_{l-1}) = \prod_{j=0}^{l-1} (a_0 + a_1 \alpha_j + \dots + a_{l-1} \alpha_{j}^{l-1}).$$

If $N(\beta_w) \models \exists a_0, \dots, a_{l-1} (P(a_0, \dots, a_{l-1}) = h_w)$, then $\ord_{\p}(w) \geq \tfrac{1-l}{l}\ord_{\p}(u)$ for all factors $\p | \mathfrak{z}_b(u)$ in $K$. Conversely if $w \in C_0(u)$ and $\ord_{\p}(w) \geq \tfrac{1-l}{l}\ord_{\p}(u)$ for all factors $\p | \mathfrak{z}_b(u)$ in $K$, then $N(\beta_w) \models \exists a_0, \dots, a_{l-1} (P(a_0, \dots, a_{l-1}) = h_w)$. \bs
\end{thm}

\begin{remark}
There will be some $w \in K_0 \setminus C_0(u)$ for which $N(\beta_w) \models \exists \overline{a}(P(\overline{a}) = h_w)$. E.g.\ when $K = \fr$ and $u \in \fr$ has non-negative degree, any $w \in \fs$ suffices. \sq
\end{remark}

For an $l$-behaved $u$, fix $\p | \mathfrak{z}_b(u)$. Define the set $\ints_l(\p,u)$ to be the subset of $K$ where if $w \in \ints_l(\p,u)$, then $\ord_{\p}(w) \geq 0$, and $u^n \in \ints_l(\p,u)$ for all $n \in \N$. 

\begin{cor}\label{bigcor}
Assume $(${\Large$\star$}$)$. The set $\ints_l(\p,u)$ is $\mathcal{L}_r$-existentially definable in $K$, with one parameter $u$. Moreover, this definition is uniform in $l$-behaved $u$.
\end{cor}

\proof
First we will rewrite ``$P(a_0, \dots, a_{l-1}) = h_w$'' as a polynomial equation over $\mathbb{F}_p(u)$  with variables in $K$. This can be done with \emph{Theorem \ref{ugly}}, as for each $w \in K$, the extension $N(\beta_w)/K$ is finite. Explicitly, we do this in three steps:
\begin{enumerate}
    \item Notice the coefficients of $P$ are in $\mathbb{F}_p(a) \subset K$, by \emph{Lemma \ref{wellillbedamned}}. In addition, $[\mathbb{F}_p(a) : \mathbb{F}_p] < \infty$ by ({\Large$\star$}).
    
    \item The extension $N(\beta_w)/N$: here we set $n = 1$, $n_2 = 0$, $n_3 = l$, $t_1 = w$, $x = \beta_w$, $g(X, w) = h_w X^l - (h_w +1)$, and $f(w, a_0,\dots, a_{l-1}) = P(a_0, \dots, a_{l-1}) - h_w$. As $u$ is by design not an $l$-th power in $K$, $h_w$ is never zero. Therefore ``$P(a_0, \dots, a_{l-1}) = h_w$'' is $\mathcal{L}_u$-existentially definable with variables in $N$.
    
    \noindent Suppose the defining polynomial is $q(w, x_1, \dots, x_k) \in \mathbb{F}_p(a)[u][\overline{X}]$.
    
    \item The extension $N/K$: here we set $n = n_2 = 0$, $n_3 = k$, $x = \delta$, $g(X) = X^l - (u+1)$, and $f(w, x_1, \dots, x_k) = q(w, x_1, \dots, x_k)$. Therefore ``$P(a_0, \dots, a_{l-1}) = h_w$'' is $\mathcal{L}_u$-existentially definable in $K$. 
\end{enumerate}

Let $\varphi(w)$ be the $\mathcal{L}_u$-existential formula given by the above process; i.e.\ $K \models \varphi(w)$ implies $\ord_{\p}(w) \geq \tfrac{1-l}{l}\ord_{\p}(u)$ for all $\p | \mathfrak{z}_b(u)$, and if $\ord_{\p}(w) \geq \tfrac{1-l}{l}\ord_{\p}(u)$ for all $\p | \mathfrak{z}_b(u)$ and $w \in C_0(u)$, then $K \models \varphi(w)$.

Now fix $\p | \mathfrak{z}_b(u)$. Let $A \defeq \{w^l u^{l-1} \mbox{ : } K \models \varphi(w)\}$. Then $x \in A$ implies $\ord_{\p}(x) \geq 0$, and if $\tfrac{x}{u^{l-1}} \in (C_0(u))^{\times l}$ and $\ord_{\p}(x) \geq 0$, then $x \in A$. Now, define:
$$\ints_l(\p,u) = \{ x_1 \cdots x_l \mbox{ : } x_1, \dots, x_l \in A \cup \mathbb{F}_p \cup \{u\}\}.$$

If $x \in \ints_l(\p,u)$ then $\ord_{\p}(x) \geq 0$ by design. By construction $u^n \in \ints_l(\p,u)$ for all $n \in \N$. Finally, $\ints_l(\p,u)$ is $\mathcal{L}_r$-existentially definable with one parameter $u$, and is uniform in $l$-behaved $u$, as desired. \bs

\subsection{$p$-behaviour}

We can deduce all of $\S \ref{integrality}$ for ``$p$-behaviour'' by modifying the extensions $N(\beta_w, \alpha)/N(\beta_w)/N/K$ to be Artin-Schreier, then modifying slightly the statements of \emph{Lemma \ref{thelem} -- Corollary \ref{bigcor}} as done by Eisentr\"{a}ger \& Shlapentokh \cite{eisenshlap}. Our underlying assumption is now the following:

\begin{customass}{({\Large$\star$}$'$)}
For some $a \in C_0$, $K$ does not contain any root of $T^p - T - a$.
\end{customass}

Take $u \in K$ $p$-behaved (recall \emph{Definition \ref{lfriend}}). In addition, fix the following notation:

\begin{itemize}
    \item Let $\delta' \in \widetilde{K}$ be a root of $T^p - T - u$.
    \item For $w \in K$, $h_w = \tfrac{w^p}{u} + \tfrac{1}{u^p}$ still.
    \item Let $\beta'_w \in \widetilde{K}$ be a root of $T^p - T - \tfrac{1}{h_w}$.
    \item Let $\alpha' \in \widetilde{K}$ be a root of $T^p - T - a$. 
\end{itemize}

Once again, consider the extensions $N(\beta'_w, \alpha')/N(\beta'_w)/N = K(\delta')/K/C(u)/\mathbb{F}_p(u)$, and $N_0(\beta'_w, \alpha')/N_0(\beta'_w)/N_0 = K_0(\delta')/K_0/C_0(u)/\mathbb{F}_p(u)$. Our main tool is the following:

\begin{lem}
\cite[Lemma 6.6]{eisenshlap}. Let $L$ be a function field of characteristic $p$. Let $z \in L$ and let $\gamma \in \widetilde{L}$ be a root of the equation $X^p - X - z = 0$. If for some $L$-prime $\mathfrak{a}$, $\ord_{\mathfrak{a}}(z) \not\equiv 0 \mod p$ and $\ord_{\mathfrak{a}}(z) < 0$, then $\mathfrak{a}$ is completely ramified in $L(\gamma)/L$. At the same time all zeros of $z$ will split completely in $L(\gamma)/L$, i.e.\ into factors of relative degree 1. \bsn
\end{lem}

Applied to the extensions in question, we recover \emph{Corollary \ref{thecor}, Lemma \ref{thelem}} (with $l = p$), \emph{Lemma \ref{thing1}} \&\ \emph{Lemma \ref{thing2}}. We deduce:

\begin{thm}
Let $\alpha_j = \alpha' +j$ for $j = 0, \dots, p-1$. Let
$$P(a_0, \dots, a_{p-1}) = \prod_{j = 0}^{p-1}(a_0 + a_1 \alpha_j + \dots + a_{p-1} \alpha_j^{p-1}).$$

If $N(\beta'_w) \models \exists \overline{a} (P(\overline{a}) = h_w)$, then $\ord_{\p}(w) \geq \tfrac{1-p}{p}\ord_{\p}(u)$ for all factors $\p | \mathfrak{z}_b(u)$ in $K$. Conversely if $w \in C_0(u)$ and $\ord_{\p}(w) \geq \tfrac{1-p}{p}\ord_{\p}(u)$ for all factors $p|\mathfrak{z}_b(u)$ in $K$, then $N(\beta'_w) \models \exists \overline{a} (P(\overline{a}) = h_w)$.
\end{thm}

\proof
Note that all solutions to an equation $X^p - X - a = 0$ in $\widetilde{K}$ can be written in the form $\alpha + i$, $i = 0, \dots, p-1$, where $\alpha$ is any solution of the equation. Also note that 
\begin{equation*}
P(a_0, \dots, a_{p-1}) = \Norm_{N(\beta'_w, \alpha_{\ast})/N(\beta'_w)}(a_0 + a_1 \alpha_{\ast} + \dots + a_{p-1} {\alpha_{\ast}}^{p-1}).
\end{equation*}
The theorem follows from the aforementioned lemmas. \bs

\begin{cor}\label{forp}
Assume $(${\Large$\star$}$')$. The set $\ints_p(\p, u)$ is $\mathcal{L}_r$-existentially definable in $K$, with one parameter $u$. Moreover, this definition is uniform in $p$-behaved $u$. \bsn
\end{cor}

\bigskip
\subsection{Assembly \&\ Initial Results.}\label{assemb}

Now we are in a position to use the above machinery. Define the unary predicate $B_{l, n}$ by $B_{l, n}(u) \Leftrightarrow u$ is $l$-behaved and for all primes $\p$, $\ord_{\p}(u) \leq n$.

\begin{thm}\label{nowthisthm}
Let $K$ be a function field of a curve, of positive characteristic and not containing the algebraic closure of a finite field. Then for some prime $l$, for $n \geq 1$ fixed, $\TThE(K; \mathcal{L}_F \cup \{B_{l, n}\})$ without parameters is undecidable.
\end{thm}

\proof
Let $C$ be the constant subfield of $K$; by assumption $C_0 = \widetilde{\mathbb{F}_p} \cap C \subsetneq \widetilde{\mathbb{F}_p}$. Hence there exists $a \in C_0$ such that either $X^p - X - a$ is irreducible, or $T^l - a$ is irreducible, for some prime $l \neq p$. If it is the former, ({\Large$\star$}$'$) is satisfied and $l = p$. If it is the latter, we may assume WLOG that $C_0$ contains an $l$-th root of unity; if not, adjoining it to $C_0$ results in a finite extension $\widehat{K}/K$. The proof below progresses for $\widehat{K}$ and we conclude the desired result for $K$ by \emph{Theorem \ref{ugly}}. Therefore ({\Large$\star$}) is satisfied with this $l$. 

If $K \models B_l(u)$ then by \emph{Corollary \ref{bigcor}/\ref{forp}} the set $\ints_l(\p, u)$ has an explicit positive-existential definition, with one parameter $u$, for a chosen fixed $\p | \mathfrak{z}_b(u)$. By \emph{Corollary \ref{pastenbuchi}}, we can give a parameter free existential $\mathcal{L}_F$-definition of the $\Den_p$ predicate. Using our above machinery, we can interpret $\ThEp(\N; 0, 1, +, |_p, <)$ in $\ThEp(K; \mathcal{L}_F)$ without parameters. Indeed, we may use the predicate $B_{l,n}$ to find some $u \in K$ such that $K \models B_{l,n}(u) \land B_{l, n}(u^n)$ -- forcing $u$ to be $l$-behaved and $\ord_{\p}(u) = 1$ for all primes $\p | \mathfrak{z}_b(u)$ (such $u \in K$ exist by Weak Approximation \cite[Theorem 2.4.1]{EP}). We then associate any natural number $a \in \N$ to the subset
$$f(a) \defeq \{w \in \ints_l(\p, u) \mbox{ : } \ord_{\p}(w) = a\},$$

where $\mbox{``}\ord_{\p}(x) = a\mbox{''} \iff \tfrac{x}{u^a}, \tfrac{u^a}{x} \in \ints_l(\p, u)$. Recall that, for all $k \in \N$, $u^k \in \ints_l(\p, u)$ -- hence $u^a \in f(a)$ for all $a \in \N$. 

The interpretation (following \cite[\S 2]{eisenshlap}) begins as follows: the equation $a = b+ c$, $a,b,c \in \N$, is equivalent to
$$\exists z_a, z_b, z_c(z_a \in f(a) \land z_b \in f(b) \land z_c \in f(c) \land z_a = z_b z_c).$$

We also have that 
\begin{align*}
a |_p b\iff \exists z_a, z_b, z ( z_a \in f(a) \land z_b \in f(b) \land \tfrac{z}{z_b}, \tfrac{z_b}{z} \in \ints_l(\p, u) \land \Den_p(z, z_a)).
\end{align*}

Finally, 
$$a \geq b \iff \exists z_a, z_b (z_a \in f(a) \land z_b \in f(b) \land \tfrac{z_a}{z_b} \in \ints_l(\p, u)).$$

This argument is uniform in $l$-behaved $u$ satisfying $K \models B_{l,n}(u) \land B_{l, n}(u^n)$. we conclude $\ThEp(K; \mathcal{L}_F \cup \{B_{l, n}\})$ without parameters is undecidable as required. \bs

\subsection{Where the method breaks down}\label{break}
One would hope that a similar undecidability result can be concluded in the language with a predicate $B_l$ for $l$-behaviour -- i.e.\ no restriction on the order of $u$. However the methods currently in existence (by which we mean interpreting arithmetic using integrality at a prime, without parameters in the style above) do not allow us to conclude undecidability via this path.

Indeed, without restriction on the order of $u$, only $\ThEp(\N; 0, +, |_p, \leq)$ is interpretable. Any nonzero constant is not interpretable (\emph{without} parameters) as a statement such as ``$\ord_{\p}(u) = k$'' is not preserved under an $\mathcal{L}_F\cup\{B_l\}$-automorphism or even monomorphism of $K$. For example, if $K = \fr$, $p$ odd, $l = 2$, then the Frobenius map given by $t \mapsto t^p$ is a $\mathcal{L}_F \cup \{B_2\}$-embedding $\fr \hookrightarrow \fr$. Therefore in general \textit{arithmetic} is not interpretable along these lines, as the constant ``1'' would be definable using multiplication.

Although $\ThEp(\N; 0, +, |_p, \leq)$ can be interpreted in $\ThE(K; \mathcal{L}_F \cup \{B_l\})$, it cannot tell us about the (un)decidability of $\ThE(K; \mathcal{L}_F \cup \{B_l\})$, as $\ThEp(\N; 0, +, |_p, \leq)$ is \textit{decidable} by an unpublished result of K.\ Kartas, communicated to the author by E.\ Hrushovski.

\bigskip
\section{Results using universal quantifiers.}\label{universal}
Some of the problems of $\S \ref{break}$ can be fixed with the introduction of an initial universal quantifier.

\begin{thm}\label{firstuni}
\textit{Let $K$ be the function field of a curve, of positive characteristic and not containing the algebraic closure of a finite field. Then for some prime $l$, $\TTh_{\forall^1\exists^+}(K; \mathcal{L}_F \cup \{\neg B_{l}\})$ without parameters is undecidable.}
\end{thm}

\proof
The proof proceeds as in \emph{Theorem \ref{nowthisthm}}, with the \emph{same interpretation} of \linebreak $\ThEp(\N; 0, 1, +, |_p, \leq)$, however now we note the constant ``1'' is only interpreted accurately when there is a prime $\p | \mathfrak{z}_b(u)$ such that $\ord_{\p}(u) = 1$.

Given a sentence $\varphi$ of $\ThEp(\N; 0, 1, +, |_p, \leq)$, let $\widehat{\varphi}(u)$ be its translation via this machinery. Consider the sentence $\forall u(B_l(u) \rightarrow \widehat{\varphi}(u))$. We claim:
$$(\N; 0, 1, +, |_p, \leq) \models \varphi \iff (K; \mathcal{L}_F \cup \{B_l\}) \models \forall u(B_l(u) \rightarrow \widehat{\varphi}(u)).$$

The reverse implication is straightforward: by Weak Approximation we may choose an $l$-behaved $u$ with $\ord_{\p}(u) = 1$ for some prime $\p | \mathfrak{z}_b(u)$. 

The forward implication requires the $\{0, 1, +, |_p, \leq\}$-embedding $\N \rightarrow \N$ given by $1 \mapsto n$ for any fixed $n \geq 1$. Therefore for any positive existential $\{0, 1, +, |_p, \leq\}$-sentence $\varphi$, written $\varphi(1)$, we have $\N \models \varphi(1) \rightarrow \forall n (n \neq 0 \rightarrow \varphi(n))$. Consequently, in the translation $K \models \forall u(B_l(u) \rightarrow \widehat{\varphi}(u))$. We conclude the desired undecidability results from the undecidability of $\ThEp(\N; 0, 1, +, |_p, \leq)$ (i.e.\ \emph{Theorem \ref{modpheidas}}). \bs

While an algebraic characterisation of $l$-behaviour seems out of reach, we are able to provide such a characterisation for ``$l$-misbehaviour''. Consider the following classical result by Leahey:

\begin{thm}\label{leahey}
\emph{\cite[Theorem, p.\ 817]{leahey}.} Let $F$ be a finite field of order $p^n$ where $p$ is a prime, $p \equiv 3 \mod 4$, and $n$ is odd. Let $f \in F[X]$ and suppose that $f = a \cdot f_1^{e_1} \cdots f_r^{e_r}$ with $a \in F$, and $f_i \in F[X]$ is the factorisation of $f$ into an element of $F$ and monic irreducible polynomials in $F[X]$.

Then $f$ can be written as the sum of two squares in $F[X]$ if and only if $e_i$ is even for those $f_i$ with odd degree. \bsn
\end{thm}

We can adapt the proof of this for the following use: let $K$ be the function field of a curve algebraic over $\mathbb{F}_p$. In this case the constant subfield $C$ is an algebraic extension of $\mathbb{F}_p$, hence $C = C_0$. We retain the assumption that $\widetilde{\mathbb{F}_p} \not\subseteq C$.

\begin{remark}\label{steinitzrem}
Note that $C$ is thus a proper subfield of $\widetilde{\mathbb{F}_p}$ and has a corresponding unique Steinitz number $\mbox{St}(C)$ (originates with Steinitz \cite{steinitz}, cf.\ \cite[Chapter 2]{inffields}). If $T^l - a$ is irreducible over $C$ for some $a \in C$ and $l \neq p$ prime, then any field extension $C'/C$ with $[C' : C] \equiv 0$ mod $l$ has solutions to $T^l = a$. Indeed, for sufficiently large $n$, $a \in \mathbb{F}_{p^n} \subset C \subset C'$, and by definition of the Steinitz number, $l | \mbox{St}(C')$ hence $\mathbb{F}_{p^{nl}} \subset C'$ and $T^l = a$ is solvable over $\mathbb{F}_{p^{nl}}$. 

A similar argument shows that if $T^p - T - a$ is irreducible over $C$ for some $a \in C$, then the polynomial $T^p - T = a$ is solvable in extensions $C'/C$ when $[C':C] \equiv 0$ \linebreak mod $p$.\sq
\end{remark}

\begin{cor}\label{usethis}
Assume $K$ satisfies either $({\Large\star})$ with $l$ prime such that $T^l - a$ is irreducible for some $a \in C$, or $({\Large\star}')$ and fix $a \in C$ such that $T^p - T - a$ is irreducible. Denote by $\alpha \in \widetilde{\mathbb{F}_p}\setminus C$ a root of the former (resp.\ the latter) polynomial.

Then $u \in K$ is not a norm  of $K(\alpha)/K$ if and only if $u$ or $\tfrac{1}{u}$ is $l$-behaved, if and only if $\tfrac{u}{u^2 + b}$ is $l$-behaved for all $b \in k\setminus\{0\}$.
\end{cor}

\proof
Suppose $u$ is $l$-behaved; it has a prime factor $\p | \mathfrak{z}(u)$ of inertial degree \&\ ramification index in $u$ both coprime to $l$. Consider the extension $K(\alpha)/K$; as it is Galois, $l = [K(\alpha):K] = e(\widehat{\p}/\p)f(\widehat{\p}/\p)g(\widehat{\p}/\p)$, where $e$ is the ramification index of $\widehat{\p} = \p \mathcal{O}_{K(\alpha)}$, $f$ is the relative inertial degree, and $g$ is the number of prime factors of $\p$ in the extension. By assumption, $[K v_{\p} : C]$ is coprime to $l$, hence the equation $X^l - a$ has no root in $K v_{\p}$. Therefore the residue field of $\widehat{\p}$ must extend that of $\p$; we conclude $e(\widehat{\p}/\p) = 1$, $f(\widehat{\p}/\p) = l$, $g(\widehat{\p}/\p) = 1$. By \emph{Lemma \ref{mention}}, as $\ord_{\p}(u)$ is coprime to $l$, $u$ is not a norm in $K(\alpha)/K$.

If $\tfrac{1}{u}$ is $l$-behaved, then we conclude by the same argument that $\tfrac{1}{u}$ is not a norm in $K(\alpha)/K$. By the multiplicative property of norms, this forces $u$ not to be a norm of $K(\alpha)/K$, as desired.

Conversely, suppose that $u$ and $\tfrac{1}{u}$ are not $l$-behaved. For all primes $\p | \mathfrak{z}(u) \cdot \mathfrak{pl}(u)$, either $\ord_{\p}(u) \equiv 0$ mod $\p$, or $[Kv_{\p} : C] \equiv 0$ mod $\p$. Let $\p_1, \dots, \p_n$ be the primes falling into the latter category. Let $\p \subset \mathcal{O}_K$ be \emph{any} prime; as the extension $K(\alpha)/K$ is unramified, either $\p$ splits completely or $f(\p\mathcal{O}_{K(\alpha)}/\p) = l$. Assume $\p = \p_i$ for some $1 \leq i \leq n$, and let $F = \mathcal{O}_K/\p = Kv_{\p}$. By \emph{Remark \ref{steinitzrem}}, $\alpha \in F$, hence $X^l \equiv a \mod \p$ is solvable over $K$, by $X = g \in \mathcal{O}_K$, say. In the extension $K(\alpha)/K$, $\p$ cannot remain irreducible, as then either $\ord_{\p}(g - \alpha) > 0$ or $\ord_{\p}(g^{l-1} + \dots + \alpha^{l-1}) > 0$ in $\mathcal{O}_K(\alpha)$. This would be a contradiction as (e.g.\ for the former) then there exists $h_0, \dots, h_{l-1} \in \mathcal{O}_K$ and $t \in \p\mathcal{O}_K$, such that $t(h_0 + \alpha h_1 + \dots + \alpha^{l-1} h_{l-1}) = g - \alpha$, hence $t h_1 = -1$; impossible as $\p$ is nontrivial. Therefore $\p$ in $K(\alpha)$ splits completely into a product of distinct primes which are conjugate over $K$. Write $\p \mathcal{O}_{K(\alpha)} = \widehat{\p}_1 \cdots \widehat{\p}_l$ and consider this locally: as 
$$l = [K(\alpha) : K] = \sum_{\q | \p} [K(\alpha)_{\q} : K_{\p}] = [K(\alpha)_{\widehat{\p}_1} : K_{\p}] + \dots + [K(\alpha)_{\widehat{\p}_l} : K_{\p}],$$
this forces $K(\alpha)_{\widehat{\p}_j} = K_{\p}$ for $1 \leq j \leq l$. 

Therefore for each $1 \leq i \leq n$, for each $1 \leq j \leq l$, $u$ is trivially a norm of $K(\alpha)_{\widehat{\p_i}_j}/K_{\p_i}$. Let $K'/\mathbb{F}_p(u)$ be a global field (with constant subfield $k'$) such that
\begin{itemize}
    \item $a\in K' \subset K$;
    \item $[K'(\alpha) : K'] = l$;
    \item In $K'$, $u$ has distinct primes $\p_1, \dots, \p_n$ such that for all $1 \leq i \leq n$, $\ord_{\p_i}(u) \not \equiv 0$ mod $l$, $[K' v_{\p_i} : k'] \equiv 0 \mod l$, and each $\p_i$ splits completely in the extension $K'(\alpha)/K'$;
    \item The primes of $\mathfrak{z}(u)$ and $\mathfrak{pl}(u)$ distinct to $\p_1, \dots \p_n$ have ramification degree divisible by $l$.
\end{itemize}

Again, this third point forces $u$ to trivially be a norm of $K'(\alpha)_{\widehat{\p_i}_j}/K'_{\p_i}$ for each $1 \leq i \leq n$, $1 \leq j \leq l$.

The extension $K'(\alpha)/K'$ is an unramified extension of global fields, hence by Weil \cite[Corollary, p.\ 226]{weil} locally every unit is a norm.  Therefore, by \emph{Lemma \ref{locallem}}, the above argument for $\p_1, \dots, \p_n$, and the Strong Hasse Principle \cite[Theorem 32.9]{reiner}, globally $u$ is a norm of $K'(\alpha)/K'$. By the same argument as previous, this forces $u$ to be a norm of $K(\alpha)/K$, as desired. 

The final equivalence is a consequence of the fact $\mathfrak{z}(\tfrac{u}{u^2+b}) = \mathfrak{z}(u) \cdot \mathfrak{pl}(u)$ for all $b \in k\setminus\{0\}$.\bs

\begin{cor}\label{firstunicor}
Let $K$ be the function field of a curve, of positive characteristic and constant subfield $C \subsetneq \widetilde{\mathbb{F}_p}$. Then $\TTh_{\forall^1\exists}(K; \mathcal{L}_F)$ without parameters is undecidable.
\end{cor}

\proof
Matching the proof of \textit{Theorem \ref{firstuni}}, we claim
$$(\N; 0, 1, +, |_p, \leq) \models \varphi \iff (K, \mathcal{L}_F) \models \forall u([\exists x \Norm_{K(\alpha)/K}(x) = u] \lor \widehat{\phi}(\tfrac{u}{u^2+1})).$$

The reverse implication is again straightforward: for some fixed prime $\p$ with $[Kv_{\p} : C] \not \equiv 0$ mod $l$, by Weak Approximation there exists an $l$-behaved $u \in K$ such that $\ord_{\p}(u) = 1$ (hence it cannot be a norm by \emph{Corollary \ref{usethis}}). So $K \models \widehat{\phi}(\tfrac{u}{u^2+1})$ where $\ord_{\p}(\tfrac{u}{u^2+1}) = 1$ still, and thus $\N \models \varphi$.

The forward implication requires the embedding argument from \textit{Theorem \ref{firstuni}}, which implied $K \models \forall u (B_l(u) \rightarrow \widehat{\phi}(u))$. In particular, $\widehat{\phi}(u)$ is satisfied by $l$-behaved elements of the form $\tfrac{u}{u^2+1}$, hence $K \models \forall u([\exists x \Norm_{K(\alpha)/K}(x) = u] \lor \widehat{\phi}(\tfrac{u}{u^2+1}))$, as required. \bs

\begin{remark}
It should be noted that this is not an entirely new result: by the work of Anscombe \& Fehm \cite{sylvy}, it is known $\Th_{\forall^1 \exists^+}(\fr; \mathcal{L}_{r})$ without parameters is undecidable (cf.\ \cite[Remark 7.9]{sylvy} where a similar argument \textit{not} using the valuation may be used). \sqn
\end{remark}

\bigskip
\subsection{Properties of $l$-behaviour.}\label{concl}

\begin{lem}\label{ptother}
Let $K$ be a function field of positive characteristic $p$ and $l$ a prime distinct to $p$. Then $u \in K$ is $l$-behaved if and only if $u^p$ is $l$-behaved, and $u \in K$ is $p$-behaved if and only if $u^l$ is $p$-behaved.
\end{lem}

\proof
Let $C$ be the field of constants of $K$. Suppose $u$ is $l$-behaved for $l \neq p$; then there exists a prime $\p$ of $K$ such that $v_{\p}(u) > 0$, $v_{\p}(u) \not \equiv 0 \mod l$ and $[K v_{\p} : C] \not\equiv 0 \mod l$. Notice $v_{\p}(u^p) = p \cdot v_{\p}(u)$, so $v_{\p}(u^p) > 0$, $v_{\p}(u^p) \not\equiv 0 \mod l$, and $[K v_{\p} : C] \not\equiv 0 \mod l$ still. Therefore $u^p$ is $l$-behaved. Conversely, if $u$ is not $l$-behaved, a similar reasoning forces us to conclude $u^p$ is not $l$-behaved.

The second statement, on $p$-behaviour, follows from the same argument as the first. This concludes the lemma. \bs

\begin{remark}\label{simvein}
Along a similar vein, notice that $u \in K$ is $l$-behaved if and only if $\sigma(u)$ is $l$-behaved, for any $\sigma \in \Aut(K)$. Indeed, this follows from three facts: first, if $\p$ is a prime of $K$, then there is a corresponding discrete valuation ring $\mathcal{O}_{\p}$ with maximal ideal $\p$. For all $\sigma \in \Aut(K)$, $\sigma(\mathcal{O}_{\p})$ is a discrete valuation ring with maximal ideal $\sigma(\p)$, hence corresponds to a prime denoted ``$\sigma(\p)$''. Second, for all $a \in K$, $\ord_{\sigma(\p)}(\sigma(a)) = \ord_{\p}(a)$. Finally, $Kv_{\p} \iso Kv_{\sigma(\p)}$. Together these facts bring us to the desired conclusion. \sq
\end{remark}

\begin{remark}\label{udisuggest}
The following is due to a suggestion of E.\ Hrushovski.

Let $k$ be a subfield of $\widetilde{\mathbb{F}_p}$, $p > 2$, such that there exists a prime $l \neq p$ with $k \setminus k^l \neq \emptyset$. Let $\mathcal{C}$ be an absolutely irreducible curve defined over $k$ of genus at least 2. We may write $k(\mathcal{C})$ as the fraction field of $k[t_1, t_2] / P$, where $P \in k[X_1, X_2]$ is an absolutely irreducible polynomial describing a smooth curve $\mathcal{C}$ of genus $\geq 2$. By assumption, for some sufficiently large finite field $\mathbb{F}_{p^n}$, $P \in \mathbb{F}_{p^n}[X_1, X_2]$.

Note that by the Riemann-Hurwitz formula, any separable endomorphism of $\mathcal{C}$ is an automorphism (any such $f : \mathcal{C} \rightarrow \mathcal{C}$ necessarily has degree 1; it follows that $f$ is injective and surjective.). Thus, any nonconstant solution of $P(X_1, X_2) = 0$ in $k(\mathcal{C})$ is of the form $(g(t_1)^{p^r}\!\!\mbox{, } g(t_2)^{p^r})$ where $g \in \Aut(K)$ and $r \geq 0$.

Let $l$ be a prime such that ({\Large$\star$}) is satisfied for $l$ and $k$ (cf.\ the beginning of \emph{Theorem \ref{nowthisthm}}). As $t_1$ is $l$-behaved in $k(\mathcal{C})$, so is $g(t_1)^{p^r}$ using \emph{Lemma \ref{ptother} \& Remark \ref{simvein}}. Therefore the first component of every nonconstant solution of $P$ in $k(\mathcal{C})$ is $l$-behaved for this $l$. The parameter-free formula $\varphi(x) = \exists v (F(x) \land F(v) \land P(x, v) = 0)$ existentially $\mathcal{L}_F$-defines a natural subset of $l$-behaved elements in $k(\mathcal{C})$. \sq
\end{remark}

\begin{lem}\label{find}
Let $K$ be a one variable algebraic function field of characteristic $p > 2$. For any nonconstant $u \in K$, there exists a prime $l \neq p$ such that $u$ is $l$-behaved. In fact, for any nonconstant $u \in K$, there exist only finitely many primes $l$ such that $u$ is \emph{not} $l$-behaved.
\end{lem}

\proof
We will show that if $[K : C(u)] \not\equiv 0 \mod l$, then $u$ is $l$-behaved. The rest of the statement of the theorem follows.

Recall $\deg(\mathfrak{z}(u)) = [K: C(u)]$. Thus, we have assumed  $\deg(\mathfrak{z}(u)) \not\equiv 0 \mod l$. In particular, there must exist a prime $\p | \mathfrak{z}(u)$ with $v_{\p}(u) > 0$, $v_{\p}(u) \not \equiv 0 \mod l$ and $[K v_{\p} : C] \not\equiv 0 \mod l$. Therefore $u$ is $l$-behaved by definition, as desired. \bs

The following is possibly known, but is easy to explicitly state as a consequence of the previously discussed machinery.

\begin{thm}\label{rephrasing}
For any prime $p > 2$, for any nonconstant $u \in \fr$, $\TThE(\mathbb{F}_{p}(t); \mathcal{L}_u)$ without parameters is undecidable, where $\mathcal{L}_u = \mathcal{L}_r \cup \{u\}$. This is to say $\TThE(\mathbb{F}_{p}(t); \mathcal{L}_r)$ with parameters in $\mathbb{F}_{p}(u)$ is undecidable.
\end{thm}

\proof
Recall that $\mathbb{F}_{p^m}$ contains a primitive $l$-th root of unity if and only if $l | p^m - 1$. We choose $m \in \N_{>0}$ large enough such that there exists a prime $l$ with $[\fr: \mathbb{F}_p(u)] \not\equiv 0 \mod l$ and $l | (p^m - 1)$. Note that necessarily $l$ is coprime to $p$, $\mathbb{F}_{p^m} \setminus (\mathbb{F}_{p^m})^l \neq \emptyset$, and by design there is a primitive $l$-th root of unity in $\mathbb{F}_{p^m}$. As $u$ does not extend the constant field of $\fr$, by Galois theory 
$$[\mathbb{F}_{p^m}(t) : \mathbb{F}_{p^m}(u)] = [\mathbb{F}_p(t) : \mathbb{F}_p(u)] \not\equiv 0 \mod l,$$

hence by \emph{Lemma \ref{find}} we conclude $u$ is $l$-behaved in the extension $\mathbb{F}_{p^m}(t)/\mathbb{F}_{p^m}(u)$. Therefore by the argument of \emph{Theorem \ref{nowthisthm}}, we deduce $\ThEp(\mathbb{F}_{p^m}(t); \mathcal{L}_u)$ without parameters is undecidable (replacing $B_{l,n}$ in the proof by $u$, and fixing $n$). As the extension $\mathbb{F}_{p^m}(t)/\fr$ is finite, $u \in \fr$, and $m$ depends solely on $u$, we may interpret $\mathbb{F}_{p^m}(t)$ in $\fr$ via \cite[\S 3]{eisenshlap} and conclude $\ThE(\mathbb{F}_{p}(t); \mathcal{L}_u)$ without parameters is undecidable, as required. \bs

\begin{cor}\label{three12}
Let $\mathcal{C}$ be a curve defined over a finite field $\mathbb{F}_p$ \emph{(}$p > 2$\emph{)} with a nonconstant $\mathbb{F}_p$-rational map $\tau : \mathcal{C} \dashrightarrow \mathbb{P}^1$. No algorithm exists which, upon input an $\mathbb{F}_p$-variety $\mathscr{A}$ and an $\mathbb{F}_p$-rational map $\rho : \mathscr{A} \dashrightarrow \mathbb{P}^1$, outputs ``YES'' if there exists an $\mathbb{F}_p$-rational map $\mathcal{C} \dashrightarrow \mathscr{A}$ making the below diagram commute, and ``NO'' otherwise. 
\[
\begin{tikzcd}
 \mathcal{C} \arrow[rr, dashed, "\exists"] \arrow[rd, dashed, "\tau", swap]& \arrow[d, phantom, "\circlearrowleft"] & \mathscr{A} \arrow[ld, dashed, "\rho"]\\
 & \mathbb{P}^1 &
\end{tikzcd}
\]
\end{cor}

\proof
If $u \in \mathbb{F}_p(t)$ is nonconstant then there exists $f \in \mathbb{F}_p[X_1, X_2]$ with $f(u, t) = 0$. We may then view $\mathbb{F}_p(t)$ as the function field over $\mathbb{F}_p$ of the curve $\mathcal{C}$ defined by $f = 0$, with subfield $\mathbb{F}_p(u)$. The result is then a geometric rephrasing of \emph{Theorem \ref{rephrasing}}, following the exposition of Poonen \cite[\S 12.1]{poonensampler}. \bsn

\vfill
\section*{Acknowledgements}
The author is grateful for the time and effort devoted to this project from Professor Ehud Hrushovski. The author would also like to extend his sincere thanks to Professor Jochen Koenigsmann, Sylvy Anscombe, Philip Dittmann, and Robert Baumann for their helpful discussions along the way. The author also extends his thanks to Professor Hector Pasten for introducing to him the paper \cite{pastennew}. Finally, the author is grateful for the terminology ``$l$-behaved'' suggested by Ronan O'Gorman.
\bigskip
\bibliography{refs}   
\bibliographystyle{acm} 
\end{document}